\documentclass[12pt,twoside]{amsart}
\usepackage{amssymb}
\usepackage{xypic} 

\title{Introduction to the theory of quasi-log varieties}
\author{Osamu Fujino} 
\date{2009/10/25, version 2.01}
\address{Department of Mathematics, Faculty of Science, Kyoto University, 
Kyoto 606-8502 Japan} 
\email{fujino@math.kyoto-u.ac.jp}
\newcommand{\Bs}[0]{{\operatorname{Bs}}}
\newcommand{\mult}[0]{{\operatorname{mult}}}
\newcommand{\Exc}[0]{{\operatorname{Exc}}}
\newcommand{\Supp}[0]{{\operatorname{Supp}}}
\newcommand{\Nqklt}[0]{{\operatorname{Nqklt}}}
\newtheorem{thm}{Theorem}[section]
\newtheorem{lem}[thm]{Lemma}
\newtheorem{cor}[thm]{Corollary}
\newtheorem{prop}[thm]{Proposition}
\newtheorem*{claim}{Claim}
\newtheorem{cla}{Claim}

\theoremstyle{definition}
\newtheorem{ex}[thm]{Example}
\newtheorem{defn}[thm]{Definition}
\newtheorem{rem}[thm]{Remark}
\newtheorem*{ack}{Acknowledgments}       
\newtheorem*{notation}{Notation}         
\newtheorem{say}[thm]{}
\begin{document}
\maketitle

\begin{abstract} 
This paper is a gentle introduction to 
the theory of quasi-log varieties by Ambro. 
We explain the fundamental theorems for the log minimal 
model program for log canonical pairs. 
More precisely, we give a proof of 
the base point free theorem for log canonical pairs in 
the framework of the theory of quasi-log varieties. 
\end{abstract}

\tableofcontents

\section{Introduction}

The aim of this article is to 
explain the fundamental 
theorems for the log minimal model program for 
log canonical pairs. 
More explicitly, we 
describe the base point free theorem 
for log canonical pairs in the framework of the theory of 
{\em{quasi-log varieties}} (see Corollary \ref{42}). 
We also treat the cone 
theorem for log canonical pairs (see Theorem \ref{cone}). 
This paper is a gentle introduction to Ambro's theory of 
quasi-log varieties (cf.~\cite{ambro}). 
It contains no new statements. 
However, it must be valuable because there are no 
introductory articles for the theory of quasi-log varieties. 
The original article \cite{ambro} seems 
to be inaccessible even for experts. 
We basically follow Ambro's arguments (see \cite
[Section 5]{ambro}) but 
we change them slightly to clarify the basic ideas 
and to remove some ambiguities and mistakes. 
The book \cite{fujino-book} is 
a comprehensive survey of the fundamental theorems 
of the log minimal model 
program from the viewpoint of the theory of 
quasi-log varieties. A new approach to the 
log minimal model program for log canonical pairs 
without 
using quasi-log varieties was found in \cite{fujino-non}. 
It seems to be more natural and much easier than the theory of 
quasi-log varieties. The paper \cite{fundamental} 
contains all the details of this new approach and is almost 
self-contained. 

Note that we only use $\mathbb Q$-divisors for simplicity. 
Some of the results can be 
generalized for $\mathbb R$-divisors with a little care. 
We do not treat the relative versions of the fundamental theorems 
in order to make our arguments transparent. 
There are no difficulties 
for the reader to obtain the relative versions 
once he understands this paper. 
We hope that this article will make the theory of quasi-log varieties 
more accessible. 
Note that 
the reader does not have to refer \cite{ambro} 
in order to read this article. 
Our formulation is 
slightly different from the one  
in \cite{ambro}. 
So, if the reader wants to taste the original 
flavor of the theory of quasi-log varieties, 
then he has to see \cite{ambro}. 

We summarize the contents of this paper. 
In Section \ref{sec2}, we quickly review the torsion-freeness and 
the vanishing theorem in \cite[Chapter 2]{fujino-book}. 
In Section \ref{sec3}, we introduce the notion of {\em{qlc pairs}}, 
which is a special case of Ambro's {\em{quasi-log varieties}}, 
and prove some important and useful lemmas. 
Theorem \ref{adj} is a key result in the theory of quasi-log varieties.  
Section \ref{sec5} is devoted to the 
proof of the base point free theorem for qlc pairs. 
This section is the heart of this paper. 
In Section \ref{sec6}, we treat the rationality 
theorem and the cone theorem for log canonical pairs. 
We note that the rationality theorem directly implies 
the essential part of the cone theorem and that 
we do not need the theory of quasi-log varieties 
for the proof of the rationality theorem. 
In the final section:~Section \ref{sec7}, we 
explain some related topics. 

\begin{ack}
I was partially supported by the Grant-in-Aid for Young Scientists 
(A) $\sharp$20684001 from JSPS. I was 
also supported by the Inamori Foundation. 
I would like to thank Takeshi Abe for 
his valuable comments. 
I would also like to thank Professor 
van der Geer and the referee for valuable suggestions and 
comments. 
\end{ack}

\subsection{Notation and Conventions}
We will work over the complex number field $\mathbb C$ throughout 
this paper. 
But we note that by using the Lefschetz principle, we can 
extend everything to the case where 
the base field is 
an algebraically closed field of characteristic 
zero. 
We will use the following notation and the notation 
in \cite{km} freely. 

\begin{notation}
(i) For a $\mathbb Q$-Weil divisor 
$D=\sum _{j=1}^r d_j D_j$ such that 
$D_j$ is a prime divisor for every $j$ and 
$D_i\ne D_j$ for $i\ne j$, we define 
the {\em{round-up}} $\ulcorner D\urcorner =\sum _{j=1}^{r} 
\ulcorner d_j\urcorner D_j$ 
(resp.~the {\em{round-down}} $\llcorner D\lrcorner 
=\sum _{j=1}^{r} \llcorner d_j \lrcorner D_j$), 
where for every rational number $x$, 
$\ulcorner x\urcorner$ (resp.~$\llcorner x\lrcorner$) is 
the integer defined by $x\leq \ulcorner x\urcorner <x+1$ 
(resp.~$x-1<\llcorner x\lrcorner \leq x$). 
The {\em{fractional part}} $\{D\}$ 
of $D$ denotes $D-\llcorner D\lrcorner$. 
We define 
$$
D^{=1}=\sum_{d_j=1} D_j, \ \ \text{and}\ \ \ 
D^{<1}=\sum_{d_j<1}d_j D_j.
$$ 
We call $D$ a {\em{boundary}} 
(resp.~{\em{subboundary}}) 
$\mathbb Q$-divisor if 
$0\leq d_j\leq 1$ (resp.~$d_j \leq 1$) for all $j$. 
Note that {\em{$\mathbb Q$-linear equivalence}} of two 
$\mathbb Q$-divisors $B_1$ and $B_2$ is denoted 
by $B_1\sim _{\mathbb Q}B_2$. 

(ii) For a proper birational morphism $f:X\to Y$, 
the {\em{exceptional locus}} $\Exc (f)\subset X$ is the locus where 
$f$ is not an isomorphism. 

(iii) Let $X$ be a normal variety and $B$ 
an effective $\mathbb Q$-divisor 
on $X$ such that $K_X+B$ is $\mathbb Q$-Cartier. 
Let $f:Y\to X$ be a resolution such that 
$\Exc (f)\cup f^{-1}_*B$ has a simple normal crossing 
support, where $f^{-1}_*B$ is the strict transform of $B$ on $Y$. 
We write $K_Y=f^*(K_X+B)+\sum _i a_i E_i$ and 
$a(E_i, X, B)=a_i$. 
We say that $(X, B)$ is {\em{lc}} if and only if 
$a_i\geq -1$ for all $i$. 
Here, lc is an abbreviation of {\em{log canonical}}. 
Note that the {\em{discrepancy}} $a(E, X, B)\in \mathbb Q$ can be 
defined for every prime divisor $E$ {\em{over}} $X$. 
Let $(X, B)$ be a lc pair. 
If $E$ is a prime divisor over $X$ such that $a(E, X, B)=-1$, 
then the center $c_X(E)$ is called a {\em{lc center}} of $(X, B)$. 
\end{notation}

\section{Vanishing and torsion-free theorems}\label{sec2}
In this section, we quickly review 
Ambro's formulation of torsion-free and vanishing 
theorems in a simplified form (see \cite[Chapter 2]{fujino-book}). 
First, we fix the notation and the conventions to state 
theorems. 

\begin{say}[Global embedded simple normal crossing 
pairs]\label{3.1}
Let $Y$ be a simple normal crossing divisor 
on a smooth 
variety $M$ and $D$ a $\mathbb Q$-divisor 
on $M$. Assume that $\Supp (D+Y)$ is simple normal crossing 
and that 
$D$ and $Y$ have no common irreducible components. 
We put $B_Y=D|_Y$ and consider the pair $(Y, B_Y)$. 
We call $(Y, B_Y)$ a {\em{global embedded simple normal 
crossing pair}}. 
Let $\nu:Y^{\nu}\to Y$ be the normalization. 
We put $K_{Y^\nu}+\Theta=\nu^*(K_Y+B_Y)$. 
A {\em{stratum}} of $(Y, B_Y)$ is an irreducible component of $Y$ or 
the image of some lc center of $(Y^\nu, \Theta^{=1})$. 
When $Y$ is smooth and $B_Y$ is a $\mathbb Q$-divisor 
on $Y$ such that 
$\Supp B_Y$ is simple normal crossing, we 
put $M=Y\times \mathbb A^1$ and $D=B_Y\times \mathbb A^1$. 
Then $(Y, B_Y)\simeq (Y\times \{0\}, B_Y\times \{0\})$ satisfies 
the above conditions, that is, 
we can consider $(Y, B_Y)$ to be a global 
embedded simple normal crossing pair. 
\end{say}
Theorem \ref{ap1} is a special 
case of the main result in \cite[Chapter 2]{fujino-book}. 
It will play crucial roles in the following 
sections. 

\begin{thm}[Torsion-freeness and vanishing 
theorem]\label{ap1}
Let $(Y, B_Y)$ be as above. 
Assume that $B_Y$ is a boundary $\mathbb Q$-divisor. 
Let 
$f:Y\to X$ be a proper morphism and $L$ a Cartier 
divisor on $Y$. 

$(1)$ Assume that $H\sim _{\mathbb Q}L-(K_Y+B_Y)$ is $f$-semi-ample.
Then, for every integer $q$,  
every non-zero local section of $R^qf_*\mathcal O_Y(L)$ contains 
in its support the $f$-image of 
some strata of $(Y, B_Y)$. 

$(2)$ Assume that $X$ is projective and 
$H\sim _{\mathbb Q}f^*H'$ for 
some ample $\mathbb Q$-Cartier 
$\mathbb Q$-divisor $H'$ on $X$. 
Then $H^p(X, R^qf_*\mathcal O_Y(L))=0$ for every $p>0$. 
\end{thm}

\begin{rem}It is obvious that 
the statement of Theorem \ref{ap1} (1) is equivalent to 
the following one. 
\begin{itemize}
\item[(1$^\prime$)] 
Assume that $H\sim _{\mathbb Q}L-(K_Y+B_Y)$ is $f$-semi-ample.
Then, for every integer $q$,  
every associated prime of $R^qf_*\mathcal O_Y(L)$ 
is the generic point of the $f$-image of 
some stratum of $(Y, B_Y)$.
\end{itemize}
\end{rem}

The above theorem follows from the next theorem. 

\begin{thm}[Injectivity theorem]\label{inj-thm} 
Let $(Y, B_Y)$ be as above.
Assume that $Y$ is proper and $B_Y$ is a boundary 
$\mathbb Q$-divisor.
Let $D$ be an effective Cartier divisor 
whose support
is contained in $\Supp \{B_Y\}$.
Assume that $L\sim _{\mathbb Q} K_Y+B_Y$.
Then the homomorphism
$$
H^q(Y, \mathcal O_Y(L))\to
H^q(Y, \mathcal O_Y(L+D)),
$$
which is induced by the natural inclusion
$\mathcal O_Y\to \mathcal O_Y(D)$, is 
injective for every $q$.
\end{thm}

For the proof, which depends on 
the theory of mixed Hodge structures, 
we recommend the reader to see \cite[Chapter 2]{fujino-book}. 
It is because \cite[Section 3]{ambro} seems to be inaccessible. 

\subsection{Idea of the proof}
We prove a very special case of Theorem \ref{inj-thm}.
This subsection is independent of the other sections.
So, the reader can skip it.
We adopt Koll\'ar's principle (cf.~\cite[Principle 2.46]{km}) 
here instead of using
the arguments by Esnault--Viehweg. 
We closely follow \cite[2.4 The Kodaira Vanishing Theorem]{km}. 
We note that \cite{fujino-jpa} 
may help the reader to understand Theorem \ref{ap1}. 
In \cite{fujino-jpa}, we give a short and 
almost self-contained proof of Theorem \ref{ap1} for the 
case when $Y$ is smooth. 

First, we recall the following Hodge theoretic results. 
Note that we compute the cohomology groups 
in the complex analytic setting throughout this subsection. 

\begin{thm}\label{hod} 
Let $V$ be a smooth projective variety and $\Sigma$ a simple
normal crossing divisor on $V$.
Let $\iota:V\setminus \Sigma\to V$ be the natural open immersion.
Then the inclusion $\iota _! \mathbb C_{V\setminus \Sigma}\subset
\mathcal O_V(-\Sigma)$ induces surjections
$$
H^i_c(V\setminus \Sigma, \mathbb C)=
H^i(V, \iota_! \mathbb C_{V\setminus \Sigma})
\to H^i(V, \mathcal O_V(-\Sigma))
$$
for all $i$.
\end{thm}
We note that $\iota_!\mathbb C_{V\setminus \Sigma}$ is 
quasi-isomorphic to 
the complex $\Omega^{\bullet}_V(\log \Sigma)\otimes
\mathcal O_V(-\Sigma)$ and the Hodge to de Rham spectral
sequence 
$$E^{p,q}_1=H^q(V, \Omega^p_V(\log \Sigma)\otimes
\mathcal O_V(-\Sigma))\Longrightarrow H^{p+q}_c(V\setminus
\Sigma, \mathbb C)
$$
degenerates at the $E_1$-term. See, for
example, \cite[I.3.]{elzein}, \cite[Section 2.4]{fujino-book}, 
or Remark \ref{2525} below. 
Theorem \ref{hod} is a direct consequence of this
$E_1$-degeneration. 

\begin{rem}\label{2525} 
We put $n=\dim V$. 
By Poincar\'e duality, 
we have $$H^{2n-(p+q)}(V\setminus \Sigma, 
\mathbb C)\simeq 
H^{p+q}_c(V\setminus \Sigma, \mathbb C)^*.$$ 
On the other hand, by Serre duality, 
we see that  
$$H^{n-q}(V, \Omega^{n-p}_V(\log \Sigma))\simeq 
H^q(V, \Omega^p_V(\log \Sigma)\otimes 
\mathcal O_V(-\Sigma))^*. $$ 
Therefore, the above $E_1$-degeneration easily 
follows from the 
well-known 
$E_1$-degeneration 
of 
$$
{}^{\prime}\!E^{n-p, n-q}_1=H^{n-q}(V, 
\Omega^{n-p}_V(\log \Sigma))\Longrightarrow 
H^{2n-(p+q)}(V\setminus \Sigma, \mathbb C). 
$$
\end{rem}

The next theorem is a special case of Theorem \ref{inj-thm} 
if we put $Y=X$, 
$L=K_X+S+M$, and $B_Y=S+\frac{d}{m}D$.

\begin{thm}\label{26-n}
Let $X$ be a smooth projective variety and $S$ 
a simple normal  crossing divisor on $X$.
Let $M$ be a Cartier divisor on $X$.
Assume that there exists a smooth divisor
$D$ on $X$ such that
$dD\sim mM$ for some relatively prime positive
integers $d$ and $m$ with $d<m$, 
$D$ and $S$ have no common irreducible 
components, and $D+S$ is a simple 
normal crossing divisor on $X$.
Then the homomorphism
$$
H^i(X, \mathcal O_X(K_X+S+M))
\to H^i(X, \mathcal O_X(K_X+S+M+bD))
$$
induced by the natural 
inclusion $\mathcal O_X\to \mathcal O_X(bD)$ is injective
for every positive integer $b$ and every $i\geq 0$.
\end{thm}

\begin{proof}
We take the usual normalization of 
the $m$-fold cyclic cover $\pi:Y\to X$
ramified along the divisor $D$ and defined by $dD\sim mM$.
We put $T=\pi^*S$.
Then $Y$ is smooth and $T$ is simple normal crossing
on $Y$.
Let $\iota:Y\setminus T\to Y$ be the natural open immersion.
Then the inclusion $\iota_!\mathbb C_{Y\setminus T}\subset 
\mathcal O_Y(-T)$ induces the following surjections
$$
H^i(Y, \iota_! \mathbb C_{Y\setminus T})\to
H^i(Y, \mathcal O_Y(-T))
$$
for all $i$ by 
Theorem \ref{hod}. 
Since the fibers of $\pi$ are zero-dimensional, 
there are no higher direct image sheaves, 
and 
$$
H^i(X, \pi_*\iota_!\mathbb C_{Y\setminus T})\to H^i(X, 
\pi_*\mathcal O_Y(-T))
$$ 
is surjective for every $i\geq 0$. 
The $\mathbb Z/m\mathbb Z$-action gives 
eigensheaf decompositions
$$
\pi_*\iota_!\mathbb C_{Y\setminus T}=
\bigoplus _{k=0}^{m-1}G_k
$$
and
$$
\pi_*\mathcal O_Y(-T)
=\bigoplus_{k=0}^{m-1} \mathcal O_X(-S-kM+\llcorner
\frac{kd}{m}\lrcorner D)
$$
such that
$$
G_k\subset \mathcal O_X(-S-kM+\llcorner
\frac{kd}{m}\lrcorner D)
$$
for $0\leq k\leq m-1$. 
By taking the $k=1$ summand, we have the surjections
$$
H^i(X, G_1)\to H^i(X, \mathcal O_X(-S-M))
$$
for all $i$.
It is easy to see that $G_1$ is a subsheaf of $\mathcal O_X(-S-M-bD)$
for every $b\geq 0$. 
See, for example, \cite[Corollary 2.54, Lemma 2.55]{km}. 
Therefore,
$$
H^i(X, \mathcal  O_X(-S-M-bD))\to 
H^i(X, \mathcal O_X(-S-M))
$$ is
surjective for every $i$ (cf.~\cite[Corollary 2.56]{km}).
By Serre duality, we have the desired injections. 
\end{proof}

By Theorem \ref{26-n}, we can easily obtain a very special 
case of Theorem \ref{ap1} (2). 
We omit the proof because it is routine work. 
See, for example, \cite[Section 2.2]{higher}. 

\begin{thm}\label{2-new} 
Let $f:X\to Y$ be a morphism from a smooth 
projective variety $X$ onto a projective variety $Y$. 
Let $S$ be a simple normal crossing divisor on $X$ and 
$L$ an ample Cartier divisor on $Y$. 
Then 
$$H^i(Y, R^jf_*\mathcal O_X(K_X+S)\otimes 
\mathcal O_Y(L))=0
$$ 
for $i>0$ and $j\geq 0$. 
\end{thm}

As a corollary, we obtain a generalization of 
the Kodaira vanishing theorem (cf.~\cite[Theorem 4.4]{fujino-jpa}). 

\begin{cor}[Kodaira vanishing theorem for log canonical 
varieties]
Let $Y$ be a projective 
variety with only log canonical singularities and 
$L$ an ample Cartier 
divisor on $Y$. 
Then 
$$
H^i(Y, \mathcal O_Y(K_Y+L))=0 
$$ 
for $i>0$. 
\end{cor}
\begin{proof}
Let $f:X\to Y$ be a resolution 
such that 
$S=\Exc (f)$ is a simple normal crossing 
divisor. 
Then $f_*\mathcal O_X(K_X+S)\simeq 
\mathcal O_Y(K_Y)$. 
Therefore, 
we have the desired 
vanishing theorem by 
Theorem \ref{2-new}. 
\end{proof}

We close this subsection with Sommese's 
example. For the details and other examples, 
see \cite[Section 2.8]{fujino-book}. 

\begin{ex}
We consider $\pi:Y=\mathbb P_{\mathbb P^1}
(\mathcal O_{\mathbb P^1}\oplus 
\mathcal O_{\mathbb P^1}(1)^{\oplus 3})\to \mathbb P^1$. 
Let $\mathcal M$ denote the tautological 
line bundle of $\pi:Y\to \mathbb P^1$. 
We take a general member 
$X$ of $|(\mathcal M\otimes 
\pi^*\mathcal O_{\mathbb P^1}(-1))^{\otimes 4}|$. 
Then $X$ is a normal Gorenstein 
projective 
threefold. 
Note that $X$ is not lc. 
We put $\mathcal O_Y(L)=\mathcal M\otimes \pi^*\mathcal O_{\mathbb P^1}
(1)$. Then $L$ is an ample Cartier divisor on $Y$. 
We can check that 
$H^1(X, \mathcal O_X(K_X+L))=\mathbb C$. Thus, 
the Kodaira vanishing theorem does not 
necessarily hold 
for non-lc varieties. 
\end{ex}

\section{Adjunction for qlc varieties}\label{sec3}

To prove the base point free theorem for 
log canonical pairs following Ambro's 
idea, 
it is better to introduce the notion of 
{\em{qlc varieties}}. For the 
details, see \cite[Section 3.2]{fujino-book}. 

\begin{defn}[Qlc varieties]\label{qlc}
A {\em{qlc variety}} is a variety $X$ with a 
$\mathbb Q$-Cartier $\mathbb Q$-divisor 
$\omega$, and a finite collection $\{C\}$ of reduced 
and irreducible subvarieties of $X$ such that there is a 
proper morphism $f:(Y, B_Y)\to X$ from a global embedded simple 
normal crossing pair as in \ref{3.1} 
satisfying the following properties: 
\begin{itemize}
\item[(1)] $f^*\omega\sim_{\mathbb Q}K_Y+B_Y$ such 
that $B_Y$ is a subboundary $\mathbb Q$-divisor. 
\item[(2)] There is an isomorphism 
$$
\mathcal O_X\simeq f_*\mathcal O_Y(\ulcorner 
-(B_Y^{<1})\urcorner). 
$$
\item[(3)] The collection of subvarieties $\{C\}$ coincides with the image 
of $(Y, B_Y)$-strata.  
\end{itemize}
We use the following terminology. 
The subvarieties $C$ 
are the {\em{qlc centers}} of $X$, 
and $f:(Y, B_Y)\to X$ is a {\em{quasi-log resolution}} 
of $X$. We sometimes simply say that 
$[X, \omega]$ is a {\em{qlc pair}}, or the 
pair $[X, \omega]$ is {\em{qlc}}. 
\end{defn}

\begin{rem}\label{32}
By condition (2), we have an isomorphism 
$\mathcal O_X\simeq f_*\mathcal O_Y$. In particular, 
$f$ is a surjective morphism with connected fibers and 
$X$ is semi-normal. 
\end{rem}

\begin{prop}\label{43} 
Let $(X, B)$ be a lc pair. 
Then $[X, K_X+B]$ is a qlc pair. 
\end{prop}
\begin{proof}
Let $f:Y\to X$ be a resolution such that 
$K_Y+B_Y=f^*(K_X+B)$ and 
$\Supp B_Y$ is a simple normal crossing 
divisor. 
Then $\mathcal O_X\simeq f_*\mathcal O_Y(\ulcorner 
-(B^{<1}_Y)\urcorner)$ because 
$\ulcorner -(B^{<1}_Y)\urcorner$ is effective and $f$-exceptional. 
We note that a qlc center $C$ is $X$ itself or a lc center 
of $(X, B)$. 
\end{proof}

We start with an easy lemma. 

\begin{lem}\label{3434}
Let $f:Z\to Y$ be a proper birational morphism between 
smooth varieties and $B_Y$ a subboundary 
$\mathbb Q$-divisor on $Y$ such 
that $\Supp B_Y$ is simple normal crossing. 
Assume that $K_Z+B_Z=f^*(K_Y+B_Y)$ and 
that $\Supp B_Z$ is simple normal crossing. 
Then we have $$f_*\mathcal O_Z(\ulcorner -(B^{<1}_Z)\urcorner 
)\simeq 
\mathcal O_Y(\ulcorner -(B^{<1}_Y)\urcorner 
). $$ 
\end{lem}

\begin{proof}
By $K_Z+B_Z=f^*(K_Y+B_Y)$, we 
obtain 
\begin{align*}
K_Z=&f^*(K_Y+B^{=1}_Y+\{B_Y\})\\&+f^*
(\llcorner B^{<1}_Y\lrcorner)
-\llcorner B^{<1}_Z\lrcorner
-B^{=1}_Z-\{B_Z\}.
\end{align*} 
If $a(\nu, Y, B^{=1}_Y+\{B_Y\})=-1$ for 
a prime divisor $\nu$ over $Y$, then 
we can check that $a(\nu, Y, B_Y)=-1$ by using 
\cite[Lemma 2.45]{km}. 
Since $f^*
(\llcorner B^{<1}_Y\lrcorner)
-\llcorner B^{<1}_Z\lrcorner$ is 
Cartier, we can easily see that 
$f^*(\llcorner B^{<1}_Y\lrcorner)
=\llcorner B^{<1}_Z\lrcorner+E$, 
where $E$ is an effective $f$-exceptional divisor. 
Thus, we obtain 
$$f_*\mathcal O_Z(\ulcorner -(B^{<1}_Z)\urcorner)\simeq 
\mathcal O_Y(\ulcorner -(B^{<1}_Y)\urcorner).$$  
This completes the proof.
\end{proof}

The following lemma is very important in the study of 
qlc pairs. 

\begin{lem}\label{3535} 
We use the same notation and assumption 
as in {\em{Lemma \ref{3434}}}. 
Let $S$ be a simple normal crossing divisor on $Y$ such 
that $S\subset \Supp B^{=1}_Y$. Let $T$ be the union of the 
irreducible 
components of $B^{=1}_Z$ that are mapped into $S$ by $f$. 
Assume that 
$\Supp f^{-1}_*B_Y\cup \Exc (f)$ is simple normal crossing on $Z$. 
Then we have 
$$f_*\mathcal O_T(\ulcorner -(B^{<1}_T)\urcorner)\simeq 
\mathcal O_S(\ulcorner -(B^{<1}_S)\urcorner 
),$$ 
where 
$(K_Z+B_Z)|_T=K_T+B_T$ and $(K_Y+B_Y)|_S=K_S+B_S$. 
\end{lem}
\begin{proof}
We use the same notation as in the proof of Lemma \ref{3434}. 
We consider the short exact sequence 
$$
0\to \mathcal O_Z(\ulcorner -(B^{<1}_Z)\urcorner-T)\to 
\mathcal O_Z(\ulcorner -(B^{<1}_Z)\urcorner)
\to \mathcal O_T(\ulcorner -(B^{<1}_T)\urcorner)\to 0. 
$$
Since $T=f^*S-F$, where $F$ is an effective $f$-exceptional 
divisor, 
we obtain 
$$
\ulcorner -(B^{<1}_Z)\urcorner-T=f^*(\ulcorner -(B^{<1}_Y)\urcorner-S)
+E+F. 
$$ 
Here, we used $f^*(\llcorner B^{<1}_Y\lrcorner)=\llcorner B^{<1}_Z\lrcorner
+E$ in the proof of Lemma \ref{3434}. 
Therefore, 
\begin{align*}
f_*\mathcal O_Z(\ulcorner -(B^{<1}_Z)\urcorner-T)&\simeq 
\mathcal O_Y(\ulcorner -(B^{<1}_Y)\urcorner-S)\otimes 
f_*\mathcal O_Z(E+F)\\ 
&\simeq \mathcal O_Y(\ulcorner -(B^{<1}_Y)\urcorner-S). 
\end{align*}
It is because $E$ and $F$ are effective and $f$-exceptional. 
We note that 
\begin{align*} 
&(\ulcorner -(B^{<1}_Z)\urcorner 
-T)-(K_Z+\{B_Z\}+(B^{=1}_Z-T))
\\ &=-f^*(K_Y+B_Y). 
\end{align*} 
Therefore, every local section 
of $R^1f_*\mathcal O_Z(\ulcorner -(B^{<1}_Z)\urcorner 
-T)$ contains in its support the $f$-image 
of some strata of $(Z, \{B_Z\}+B^{=1}_Z-T)$ by 
Theorem \ref{ap1} (1). 
\begin{claim}
No strata of $(Z, \{B_Z\}+B^{=1}_Z-T)$ are 
mapped into $S$ by $f$. 
\end{claim}
\begin{proof}[Proof of Claim]
Assume that there is a stratum $C$ of $(Z, \{B_Z\}+B^{=1}_Z-T)$ such that 
$f(C)\subset S$. Note that 
$\Supp f^*S\subset \Supp f^{-1}_*B_Y\cup \Exc (f)$ and 
$\Supp B^{=1}_Z\subset \Supp f^{-1}_*B_Y\cup \Exc (f)$. 
Since $C$ is also a stratum of $(Z, B^{=1}_Z)$ and 
$C\subset \Supp f^*S$, 
there exists an irreducible component $G$ of $B^{=1}_Z$ such that 
$C\subset G\subset \Supp f^*S$. 
Therefore, by the definition of $T$, $G$ is an 
irreducible component of $T$ because $f(G)\subset S$ and $G$ is an 
irreducible component of $B^{=1}_Z$. So, 
$C$ is not a stratum of $(Z, \{B_Z\}+B^{=1}_Z-T)$. It is 
a contradiction. 
\end{proof}
On the other hand, $f(T)\subset S$. Therefore, 
$$f_*\mathcal O_T(\ulcorner -(B^{<1}_T)\urcorner
)
\to R^1f_*\mathcal O_Z(\ulcorner -(B^{<1}_Z)\urcorner
-T)$$ is a zero map by 
the above claim. 
Thus, we obtain 
$$f_*\mathcal O_T(\ulcorner -(B^{<1}_T)\urcorner 
)\simeq 
\mathcal O_S(\ulcorner -(B^{<1}_S)\urcorner 
)$$ 
by the following commutative diagram. 
$$
\xymatrix{
0\ar[r]
&\mathcal O_Y(\ulcorner -(B^{<1}_Y)\urcorner-S)
\ar[d]^{\simeq}\ar[r]
& \mathcal O_Y(\ulcorner -(B^{<1}_Y)\urcorner)
\ar[d]^{\simeq}\ar[r]
&\mathcal O_S(\ulcorner -(B^{<1}_S)\urcorner)
\ar[d]\ar[r]&0\\
0\ar[r]
&f_*\mathcal O_Z(\ulcorner -(B^{<1}_Z)\urcorner-T)
\ar[r]&
f_*\mathcal O_Z(\ulcorner -(B^{<1}_Z)\urcorner)
\ar[r]
&f_*\mathcal O_T(\ulcorner -(B^{<1}_T)\urcorner)
\ar[r]&0
}
$$ 
This completes the proof. 
\end{proof}

The following theorem (cf.~\cite[Theorem 4.4]{ambro}) 
is one of the key results 
for the theory of qlc varieties. 
It is a consequence of Theorem \ref{ap1}. 
See also Theorem \ref{thm666} below. 

\begin{thm}[Adjunction and vanishing theorem]\label{adj}
Let $[X, \omega]$ be a qlc pair and $X'$ 
a union of some qlc centers of $[X, \omega]$. 
\begin{itemize}
\item[(i)] Then 
$[X', \omega']$ is a qlc pair, where $\omega'=\omega|_{X'}$. 
Moreover, the qlc centers of $[X', \omega']$ 
are exactly the qlc centers of 
$[X, \omega]$ that are included in $X'$. 
\item[(ii)] Assume that $X$ is projective. 
Let $L$ be a Cartier divisor on $X$ such that 
$L-\omega$ is ample. 
Then $H^q(X, \mathcal O_X(L))=0$ and $H^q(X, \mathcal I_{X'}\otimes 
\mathcal O_X(L))=0$ for $q>0$, where $\mathcal I_{X'}$ is 
the defining ideal sheaf of $X'$ on $X$. 
Note that $H^q(X', \mathcal O_{X'}(L))=0$ for every $q>0$ 
because $[X', \omega']$ is a qlc pair by {\em{(i)}} and $L|_{X'}-\omega'$ is 
ample. 
\end{itemize}
\end{thm}

\begin{proof}
(i) Let $f:(Y, B_Y)\to X$ be a quasi-log resolution. 
Let $M$ be the ambient space 
of $Y$ and $D$ a subboundary $\mathbb Q$-divisor 
on $M$ such that $B_Y=D|_Y$. 
By taking blow-ups of $M$, 
we can assume that 
the union of all strata of $(Y, B_Y)$ mapped into 
$X'$, which is denoted by $Y'$, is a union 
of irreducible 
components of $Y$ (cf.~Lemma \ref{3535}). We put 
$Y''=Y-Y'$. 
We define $(K_Y+B_Y)|_{Y'}=K_{Y'}+B_{Y'}$ and 
consider $f:(Y', B_{Y'})\to X'$. 
We claim that $[X', \omega']$ is a qlc pair, where $\omega'=\omega|_{X'}$, 
and 
$f:(Y', B_{Y'})\to X'$ is a quasi-log resolution. 
By the definition, $B_{Y'}$ is a subboundary and 
$f^*\omega'\sim _{\mathbb Q}K_{Y'}+B_{Y'}$ on $Y'$. 
We consider the following short exact sequence 
$$
0\to \mathcal O_{Y''}(-Y')\to \mathcal O_Y\to \mathcal O_{Y'}\to 0. 
$$ 
We put $A=\ulcorner -(B^{<1}_Y)\urcorner$. 
Then we have 
$$
0\to \mathcal O_{Y''}(A-Y')\to \mathcal O_Y(A)\to \mathcal O_{Y'}(A)\to 0. 
$$ 
Applying $f_*$, we obtain 
\begin{align*}
0\to f_*\mathcal O_{Y''}(A-Y')\to \mathcal O_X\to 
f_*\mathcal O_{Y'}(A)\to R^1f_*\mathcal O_{Y''}(A-Y')\to \cdots. 
\end{align*}
The support of every non-zero local section of $R^1f_*\mathcal O_{Y''}(A
-Y')$ can not be contained in $f(Y')=X'$ by Theorem \ref{ap1} (1). 
We note that 
\begin{align*}
-f^*\omega\sim _{\mathbb Q}(A-Y')|_{Y''}-
(K_{Y''}+\{B_{Y''}\}+B^{=1}_{Y''}-Y'|_{Y''})
\end{align*} on 
$Y''$, where 
$(K_Y+B_Y)|_{Y''}=K_{Y''}+B_{Y''}$, 
and that $Y'|_{Y''}$ is contained in $B^{=1}_{Y''}$. 
Therefore, $f_*\mathcal O_{Y'}(A)\to R^1f_*\mathcal O_{Y''}(A-Y')$ is 
a zero map. 
We note that the surjection $\mathcal O_X\to f_*\mathcal O_{Y'}(A)$ decomposes 
as 
$$
\mathcal O_X\to \mathcal O_{X'}\to f_*\mathcal O_{Y'}\to f_*\mathcal 
O_{Y'}(A)
$$ 
since $f(Y')=X'$. 
Therefore, we obtain 
$$\mathcal O_{X'}\simeq f_*\mathcal O_{Y'}(A)= 
f_*\mathcal O_{Y'}(\ulcorner -(B_{Y'}^{<1})\urcorner). $$ 
Thus, we see that 
$f_*\mathcal O_{Y''}(A-Y')\simeq \mathcal I_{X'}$, the defining 
ideal sheaf of $X'$ on $X$. 
The statement for qlc centers is obvious 
by the construction of the quasi-log resolution. 
So, we obtain (i). 

(ii) Let $f:(Y, B_Y)\to X$ be a quasi-log resolution as in the proof of 
(i). Apply Theorem \ref{ap1} (2). Then 
we obtain $H^q(X, \mathcal O_X(L))=0$ for 
every $q>0$ because 
\begin{align*}
f^*(L-\omega)&\sim_{\mathbb Q}f^*L-(K_Y+B_Y)\\&=
f^*L+\ulcorner -(B^{<1}_Y)\urcorner-(K_Y+\{B_Y\}
+B^{=1}_Y)
\end{align*} 
and $f_*\mathcal O_Y(f^*L+\ulcorner -(B^{<1}_Y)
\urcorner)\simeq \mathcal O_X(L)$. 
We consider $f:Y''\to X$. 
We put $(K_Y+B_Y)|_{Y''}=K_{Y''}+B_{Y''}$. 
Then 
\begin{align*}
f^*(L-\omega)&\sim _{\mathbb Q} 
(f^*L-(K_Y+B_Y))|_{Y''}
\\ &=(f^*L+A-Y')|_{Y''}-(K_{Y''}+\{B_{Y''}\}+B^{=1}_{Y''}-Y'|_{Y''})
\end{align*} on $Y''$. 
Note that $Y'|_{Y''}$ is contained in $B^{=1}_{Y''}$. 
Therefore, 
we obtain 
$$H^q(X, f_*\mathcal O_{Y''}(A-Y')\otimes 
\mathcal O_X(L))=0$$ for every $q>0$ by Theorem \ref{ap1} (2). 
Thus this completes the proof 
by $f_*\mathcal O_{Y''}(A-Y')\simeq 
\mathcal I_{X'}$ obtained in the proof of (i). 
\end{proof}

\begin{cor}\label{46}
Let $[X, \omega]$ be a qlc pair and 
$X'$ an irreducible component of $X$. 
Then $[X', \omega']$, where $\omega'=\omega|_{X'}$, 
is a qlc pair. 
\end{cor}
\begin{proof}
By Definition \ref{qlc} and Remark \ref{32}, 
$X'$ is a qlc center of $[X, \omega]$. 
Therefore, by Theorem \ref{adj} (i), $[X', \omega']$ 
is a qlc pair. 
\end{proof}

We use the next definition in Section \ref{sec5}. 

\begin{defn}\label{3838}
Let $[X, \omega]$ be a qlc pair. 
Let $X'$ be the union of qlc centers of $X$ that are not 
any irreducible components of $X$. 
Then $X'$ with $\omega'=\omega|_{X'}$ is a qlc variety by Theorem 
\ref{adj} (i). 
We denote it by $\Nqklt (X, \omega)$. 
\end{defn}

We close this section with the following very 
useful 
lemma, which seems to be indispensable 
for the proof of the base point free theorem in Section \ref{sec5}. 

\begin{lem}\label{47} 
Let $f:(Y, B_Y)\to X$ be a quasi-log resolution of a qlc pair 
$[X, \omega]$. 
Let $E$ be a Cartier divisor on $X$ such that 
$\Supp E$ contains no qlc centers 
of $[X, \omega]$. 
By blowing up $M$, the ambient space of 
$Y$, 
inside $\Supp f^*E$, 
we can assume that $(Y, B_Y+f^*E)$ is a global 
embedded simple normal crossing 
pair. 
\end{lem}
\begin{proof}
First, we take a blow-up of $M$ along $f^*E$ and apply Hironaka's 
resolution theorem to $M$. 
Then we can assume that there exists a Cartier divisor 
$F$ on $M$ such that 
$\Supp (F\cap Y)=\Supp f^*E$. 
Next, we apply Szab\'o's resolution lemma 
to $\Supp (D+Y+F)$ on $M$. 
Thus, we obtain the desired 
properties by Lemma \ref{3535}. 
\end{proof}

\section{Base point free theorem}\label{sec5}
The next theorem is the main theorem of this section. 
It is a special case of \cite[Theorem 5.1]{ambro}. 
This formulation is indispensable for 
the inductive treatment of log canonical pairs 
in the framework of the theory of 
quasi-log varieties. 
For the details, see \cite[Section 3.2.2]{fujino-book}. 

\begin{thm}\label{bpfree} 
Let $[X, \omega]$ be a projective qlc pair and 
$L$ a nef Cartier divisor on $X$. 
Assume that $qL-\omega$ is ample for some 
$q>0$. 
Then $\mathcal O_X(mL)$ is generated by global 
sections for every $m\gg0$, 
that is, there exists a positive number 
$m_0$ such that 
$\mathcal O_X(mL)$ is generated by global sections for 
every $m\geq m_0$.  
\end{thm}

\begin{proof} 
First, we note that the statement is obvious when $\dim X=0$. 

\begin{cla}\label{c1}
We can assume that $X$ is irreducible. 
\end{cla}
Let $X'$ be an irreducible component of $X$. 
Then $X'$ with $\omega'=\omega|_{X'}$ 
has a natural qlc structure induced by $[X, \omega]$ by adjunction (see 
Corollary \ref{46}). 
By the vanishing theorem (see Theorem \ref{adj} (ii)), we have 
$H^1(X, \mathcal I_{X'}\otimes 
\mathcal O_X(mL))=0$ for all $m\geq q$. 
We consider the following 
commutative diagram. 
$$
\xymatrix{
H^0(X, \mathcal O_X(mL))\otimes \mathcal O_X 
\ar[d]\ar[r]^{\alpha}& 
H^0(X', \mathcal O_{X'}(mL))\otimes \mathcal O_{X'}
\ar[d]\ar[r]&0\\
\mathcal O_X(mL)
\ar[r]&
\mathcal O_{X'}(mL)
\ar[r]&0
}
$$
Since $\alpha$ is surjective for $m\geq q$, we can assume that 
$X$ is irreducible when we prove this theorem. 

\begin{cla}\label{c2} 
For every $m\gg0$, 
$\mathcal O_X(mL)$ is generated by global sections 
on an open neighborhood of $\Nqklt(X, \omega)$. 
\end{cla}

We put $X'=\Nqklt(X, \omega)$. 
Then $[X', \omega']$, 
where $\omega'=\omega|_{X'}$, is a qlc pair 
by adjunction (see Definition \ref{3838} and 
Theorem \ref{adj} (i)). 
By the induction on the dimension, 
$\mathcal O_{X'}(mL)$ 
is generated by global sections for every $m\gg 0$. 
By the following commutative diagram: 
$$
\xymatrix{
H^0(X, \mathcal O_X(mL))\otimes \mathcal O_X 
\ar[d]\ar[r]^{\alpha}& 
H^0(X', \mathcal O_{X'}(mL))\otimes \mathcal O_{X'}
\ar[d]\ar[r]&0\\
\mathcal O_X(mL)
\ar[r]&
\mathcal O_{X'}(mL)
\ar[r]&0, 
}
$$
we know that, for every $m\gg0$, 
$\mathcal O_X(mL)$ is generated 
by global sections on an open neighborhood of $X'$. 

\begin{cla}\label{c3} 
For every $m\gg0$, 
$\mathcal O_X(mL)$ is generated by global sections 
on a non-empty Zariski open set. 
\end{cla} 

By Claim \ref{c2}, we can assume that $\Nqklt(X, \omega)$ is empty. 
If $L$ is numerically trivial, then 
$H^0(X, \mathcal O_X(L))=H^0(X, \mathcal O_X(-L))=\mathbb C$. 
It is because $h^0(X, \mathcal O_X(\pm L))=\chi (X, \mathcal O_X(\pm L))
=\chi (X, \mathcal O_X)=1$  
by Theorem \ref{adj} (ii) and 
\cite[Chapter II \S2 Theorem 1]{kleiman}. 
Therefore, $\mathcal O_X(L)$ is trivial. 
So, we can assume that $L$ is not numerically trivial. 
Let $f:(Y, B_Y)\to X$ be a quasi-log resolution. 
Let $x\in X$ be a general smooth point. Then 
we can take a $\mathbb Q$-divisor $D$ such that 
$\mult _x D> \dim X$ and that 
$D\sim _{\mathbb Q}(q+r)L-\omega$ for some $r>0$ 
(see \cite[3.5 Step 2]{km}). 
By blowing up $M$, 
we can assume that $(Y, B_Y+f^*D)$ is a global embedded 
simple normal crossing pair by Lemma \ref{47}. 
We note that every stratum of $(Y, B_Y)$ is mapped onto $X$ by the 
assumption. 
By the construction of $D$, 
we can find a positive rational number $c<1$ 
such that $B_Y+cf^*D$ is a subboundary 
and some stratum of $(Y, B_Y+cf^*D)$ does not dominate $X$. 
Here, we used $f_*\mathcal O_Y(\ulcorner 
-(B^{<1}_Y)\urcorner)\simeq \mathcal O_X$. 
Then the pair $[X, \omega+cD]$ is qlc and $f:(Y, B_Y+cf^*D)\to 
X$ is a quasi-log resolution. 
We note that $q'L-(\omega+cD)$ is ample by $c<1$, where 
$q'=q+cr$. 
By the construction, $\Nqklt(X, \omega+cD)$ is non-empty. 
Therefore, by applying Claim \ref{c2} to 
$[X, \omega+cD]$, for every $m\gg0$, 
$\mathcal O_X(mL)$ is 
generated by 
global sections on an open neighborhood of 
$\Nqklt(X, \omega+cD)$. 
So, we obtain Claim \ref{c3}. 

Let $p$ be a prime number and $l$ a large integer. 
Then $|p^lL|\ne \emptyset$ by Claim \ref{c3} and 
$|p^l L|$ is free on an open neighborhood of 
$\Nqklt(X, \omega)$ by Claim \ref{c2}. 

\begin{cla}\label{c4} 
If the base locus $\Bs|p^lL|$ $($with reduced scheme 
structure$)$ is not empty, 
then $\Bs|p^{l'}L|$ is strictly smaller than $\Bs|p^{l}L|$ for 
some $l'> l$. 
\end{cla}

Let $f:(Y, B_Y)\to X$ be a quasi-log resolution. 
We take a general member $D\in |p^lL|$. We note 
that $|p^lL|$ is free on an open neighborhood of 
$\Nqklt(X, \omega)$. 
Thus, 
$f^*D$ intersects all strata of $(Y, \Supp B_Y)$ transversally 
over $X\setminus \Bs|p^lL|$ by Bertini and 
$f^*D$ contains no strata of $(Y, B_Y)$. 
By taking blow-ups of $M$ suitably, 
we can assume that $(Y, B_Y+f^*D)$ is a global 
embedded simple normal crossing pair (cf.~Lemmas \ref{47} 
and \ref{3535}). 
We take the maximal positive rational number $c$ such 
that $B_Y+cf^*D$ is a subboundary. 
We note that $c\leq 1$. 
Here, we used $\mathcal O_X\simeq 
f_*\mathcal O_Y(\ulcorner -(B^{<1}_Y)\urcorner)$. 
Then $f:(Y, B_Y+cf^*D)\to X$ is a quasi-log resolution 
of $[X, \omega'=\omega+cD]$. 
Note that $[X, \omega']$ has a qlc center $C$ that 
intersects $\Bs|p^lL|$ by the construction. 
By the induction 
on the dimension, $\mathcal O_C(mL)$ is generated by global sections 
for all $m\gg 0$. 
We can lift the sections of $\mathcal O_C(mL)$ to $X$ for 
$m\geq q+cp^l$ by Theorem \ref{adj} (ii). Then 
we obtain that, for every $m\gg0$, 
$\mathcal O_X(mL)$ is generated by 
global sections on an open neighborhood of $C$. 
Therefore, $\Bs|p^{l'}L|$ is strictly smaller 
than $\Bs|p^{l}L|$ for some $l'> l$. 

\begin{cla}\label{c5} 
$\mathcal O_X(mL)$ is generated by global sections 
for every $m\gg 0$. 
\end{cla}
By Claim \ref{c4} and 
the noetherian induction, 
$\mathcal O_X(p^{l}L)$ and $\mathcal O_X({p'}^{l'}L)$ 
are generated by global sections for large $l$ and $l'$, 
where $p$ and $p'$ are prime numbers and $p\ne p'$. 
So, there exists a positive number 
$m_0$ such that $\mathcal O_X(mL)$ is 
generated by global sections 
for every $m\geq m_0$. 
\end{proof}

The next corollary is obvious by Theorem \ref{bpfree} 
and Proposition \ref{43}. 

\begin{cor}[Base point free theorem for lc pairs]\label{42} 
Let $(X, B)$ be a projective lc pair and $L$ 
a nef Cartier divisor on $X$. 
Assume that $qL-(K_X+B)$ is ample 
for some $q>0$. 
Then $\mathcal O_X(mL)$ is generated by global sections for every 
$m\gg 0$.  
\end{cor}

The reader can find another proof of Corollary \ref{42} 
in \cite[Section 4]{fujino-non}. 
It does not need the notion of qlc pairs. 

\section{Cone theorem}\label{sec6}

In this section, we will state the cone theorem 
for lc pairs (cf.~Theorem \ref{cone}). The 
essential part of the cone theorem follows from the rationality 
theorem:~Theorem \ref{61}. 
The rationality theorem is in turn implied by the vanishing theorem 
for lc centers (cf.~Theorem \ref{thm666}) by 
the standard argument (for the details, see \cite[Section 5]{fujino-non}). 
Note that Theorem \ref{thm666} is a special case of Theorem \ref{adj} (ii), but it can be proved much more easily 
(see, for example, \cite[Theorem 4.1]{fujino-jpa} or 
\cite[Theorem 2.2]{fujino-non}). 
Note that we do not need the theory of quasi-log varieties 
in this section. 
So, we omit the details.

\subsection{Rationality theorem} 
Here, we explain the rationality theorem 
for log canonical pairs. 
It implies the essential part of 
the cone theorem for log canonical 
pairs. 

\begin{thm}[Rationality theorem]\label{61} 
Let $(X, B)$ be a projective lc pair such that 
$a(K_X+B)$ is Cartier for a positive integer $a$. 
Let $H$ be an ample Cartier divisor 
on $X$. Assume that $K_X+B$ is not nef. 
We put 
$$
r=\max \{ t \in \mathbb R : H+t(K_X+B)  \text{ is nef}\,  \}. 
$$ 
Then $r$ is a rational number of the 
form $u/v$ $($$u, v\in \mathbb Z$$)$ where 
$0<v\leq a(\dim X+1)$. 
\end{thm}

As we explained above, Theorem \ref{61} can be proved 
easily by using the following very special case of 
Theorem \ref{adj} (ii). 

\begin{thm}[Vanishing theorem for lc centers]\label{thm666} 
Let $X$ be a projective variety and 
$B$ a boundary $\mathbb Q$-divisor 
on $X$ such that $(X, B)$ is log canonical. 
Let $D$ be a Cartier divisor on $X$. Assume that 
$D-(K_X+B)$ is ample. 
Let $C$ be a lc center of the pair $(X, B)$ with a 
reduced scheme structure. 
Then we have 
\begin{align*}
H^i(X, \mathcal I_C\otimes \mathcal O_X(D))=0, \ \ 
H^i(C, \mathcal O_C(D))=0
\end{align*}
for all $i>0$, where 
$\mathcal I_C$ is the defining ideal sheaf of 
$C$ on $X$. In particular, 
the restriction map 
$$
H^0(X, \mathcal O_X(D))\to H^0(C, \mathcal O_C(D))
$$  
is surjective. 
\end{thm}

The reader can find the details of the rationality theorem 
in \cite[Section 5]{fujino-non}. 

\subsection{Cone theorem}

Let us state the main theorem of this section. 

\begin{thm}[Cone theorem]\label{cone} 
Let $(X, B)$ be a projective 
lc pair. 
Then we have 
\begin{itemize}
\item[(i)] There are $($countably many$)$ 
rational curves $C_j\subset X$ such that $0<-(K_X+B)\cdot C_j\leq 
2\dim X$, and 
$$
\overline {NE}(X)=\overline {NE}(X)_{(K_X+B)\geq 0}
+\sum \mathbb R_{\geq 0}[C_j]. 
$$
\item[(ii)] For any $\varepsilon >0$ and ample 
$\mathbb Q$-divisor $H$, 
$$
\overline {NE}(X)=\overline {NE}(X)_{(K_X+B+\varepsilon 
H)\geq 0}
+\underset{\text{finite}}\sum \mathbb R_{\geq 0}[C_j]. 
$$
\item[(iii)] Let $F\subset \overline {NE}(X)$ be a 
$(K_X+B)$-negative 
extremal face. Then there is a unique morphism 
$\varphi_F:X\to Z$ such that 
$(\varphi_F)_*\mathcal O_X\simeq \mathcal O_Z$, 
$Z$ is projective, and 
an irreducible curve $C\subset X$ is mapped to a point 
by $\varphi_F$ if and only if $[C]\in F$. 
The map 
$\varphi_F$ is called the {\em{contraction}} of $F$. 
\item[(iv)] Let $F$ and $\varphi_F$ be as in {\em{(iii)}}. Let $L$ 
be a line bundle on $X$ such that 
$(L\cdot C)=0$ for every curve $C$ with 
$[C]\in F$. Then there is a line bundle $L_Z$ on $Z$ such that 
$L\simeq \varphi_F^*L_Z$. 
\end{itemize}
\end{thm}

\begin{proof}
The estimate $\leq 2\dim X$ and 
the fact that $C_j$ is a rational curve in (i) can be proved 
by Kawamata's argument in \cite{kawamata} with the aid of \cite{bchm}. 
For the details, see \cite[Section 3.1.3]{fujino-book} or 
\cite[Section 18]{fundamental}. 
The other statements in (i) and (ii) are formal consequences of 
the rationality theorem (cf.~Theorem \ref{61}). For the proof, 
see \cite[Theorem 3.15]{km}. 
The statements (iii) and (iv) are obvious by Corollary \ref{42} and 
the statements (i) and (ii). See Steps 7 and 9 in \cite[3.3 The Cone 
Theorem]{km}. 
\end{proof}

\section{Related topics}\label{sec7}  

In this paper, we did not prove Theorem \ref{ap1}, which 
is a key result for the theory of quasi-log varieties. 
For the proof, see \cite[Chapter 2]{fujino-book}. 
The paper \cite{fujino-jpa} is a gentle introduction 
to the vanishing and torsion-free theorems. 
In \cite[Chapters 3, 4]{fujino-book}, we gave a proof of the 
existence of fourfold lc flips and 
proved the base point free theorem of Reid--Fukuda type 
for lc pairs. 
The base point free theorem for 
lc pairs was generalized in \cite{fujino3}, where 
we obtained Koll\'ar's effective base point free theorem 
for lc pairs. 
In \cite{fujino4}, we proved 
the effective base point free theorem of Angehrn--Siu 
type for lc pairs. 
Recently, we introduced the notion of non-lc ideal 
sheaves and proved the restriction theorem (see \cite{fujino5}). 
It is a generalization of 
Kawakita's inversion of adjunction on log canonicity for 
normal divisors. In \cite{fujino-nagata}, 
we proved that the log canonical ring is finitely generated in 
dimension four. 
In \cite{fujino-non}, we succeeded in proving 
the fundamental theorems of 
the log minimal model program for log canonical 
pairs without using the theory of quasi-log varieties. 
Our new approach in \cite{fujino-non} seems to 
be more natural and simpler than Ambro's theory of 
quasi-log varieties. 
In \cite{fundamental}, 
we went ahead with this new approach. We strongly recommend the 
reader to see \cite{fujino-non} and 
\cite{fundamental}. 

\ifx\undefined\bysame
\newcommand{\bysame|{leavemode\hbox to3em{\hrulefill}\,}
\fi

\end{document}